\documentclass [11pt,twoside,a4paper]{article}
\usepackage{amsfonts}
\usepackage{amsthm}
\usepackage{amsmath}
\usepackage{amstext}
\usepackage{amssymb}
\usepackage{mathrsfs}
\usepackage{amscd}
\usepackage{xypic}
\usepackage{epsf}              
\usepackage{graphicx}          
\usepackage{fancybox}          
\usepackage{color}             
\usepackage{fancyhdr}
\usepackage[hang,footnotesize]{caption2}  

\def\mathcal{\mathscr}
\newfont{\aaa}{cmb10 at 16pt}
\newfont{\bbb}{cmb10 at 11pt}
\newtheorem{lem}{Lemma}[]
\newtheorem{thm}{Theorem}
\newtheorem{cor}{Corollary}[]
\newtheorem{rem}{Remark}[]

\newtheorem{defn}{Definition}[]
\def\v1{\vspace{1mm}}

\def\leq{\leqslant}

\def\geq{\geqslant}
\def\N{\Bbb N}

\def\R{{\mathbb R}}

\def\E{{\mathbb E}}
\def\P{{\mathbb P}}

\newcommand{\beqlb}{\begin{eqnarray}}
\newcommand{\eeqlb}{\end{eqnarray}}
\newcommand{\beqnn}{\begin{eqnarray*}}
\newcommand{\eeqnn}{\end{eqnarray*}}
\newcommand{\nn}{\nonumber}

\makeatletter
\def\@evenhead{
\vbox{\hbox to \textwidth {}{\hspace{0mm}{\footnotesize
\thepage}}{\hspace{9cm} {\footnotesize {Hongshuai DAI et al.}}} \protect\vspace{1truemm}\relax \hrule depth0pt
height0.15truemm width\textwidth}}
\def\@evenfoot{}
\def\@oddhead{\vbox{\hbox to \textwidth
{{\hspace{0cm}{\footnotesize Limit theorems for functionals of Gaussian vectors}\hfill{\footnotesize
\thepage}}\hspace{0mm}}{} \protect\vspace{1truemm}\relax\hrule
depth0pt height0.15truemm width\textwidth}}
\def\@oddfoot{}
\makeatother

\begin{document}

\thispagestyle{empty} \thispagestyle{fancy} {
\fancyhead[lO,RE]{\footnotesize  XXXX \\
XXXX\\[3mm]
\includegraphics[0,-50][0,0]{11.bmp}}
\fancyhead[RO,LE]{\scriptsize \bf 
} \fancyfoot[CE,CO]{}}
\renewcommand{\headrulewidth}{0pt}


\setcounter{page}{1}
\qquad\\[8mm]

\noindent{\aaa{Limit theorems for functionals of Gaussian vectors}}\\[1mm]

\noindent{\bbb Hongshuai DAI$^{1}$,\quad Guangjun SHEN$^{2}$,\quad Lingtao KONG$^{1}$ }\\[-1mm]

\noindent\footnotesize{$1 $\quad School of Statistics, Shandong University of Finance and Economics, Jinan 250014,  China.
\\$2$\quad  Department of Mathematics, Anhui Normal  University,
 Wuhu 241000, China.}\\[6mm]

\normalsize\noindent{\bbb Abstract}\quad Operator self-similar processes, as an extension of  self-similar processes, have been studied extensively.  In this work, we study limit theorems for functionals of  Gaussian  vectors. Under some conditions, we determine that the limit of partial sums of functionals of a stationary Gaussian  sequence of random vectors  is an operator  self-similar  process.\vspace{0.3cm}

\noindent{\bbb Keywords}\quad Gaussian vector, operator self-similar process, operator fractional Brownian motion,  scaling limit\\
{\bbb MSC}\quad 60G15, 60F17\\[0.4cm]

\noindent{\bbb 1\quad Introduction}\\[0.1cm]
Self-similar processes, first studied rigorously by Lamperti \cite{LL} under the name
``semi-stable'', are stochastic processes that are invariant in distribution under suitable
scaling of time and space.  We refer to Vervaat
\cite{VW} for general properties, to Samorodnitsky and Taqqu
\cite[Chaps.7 and 8]{ST94} for studies on Gaussian and stable
self-similar processes and random fields.
 Scholars have  extended
the definition of self-similarity to allow for scaling by linear operators on $\R^d$. Let $End (\R^{d}) $ be the set of linear operators on $\R^d$
(endomorphisms) and  $Aut (\R^d)$ be the set of invertible linear
operators (automorphisms) in $End (\R^d)$. For convenience, we do
not distinguish an operator $D\in End (\R^d)$ from its associated
matrix relative to the standard basis of $\R^d$. Recall that an $\R^d$-valued stochastic  process $\tilde{Y}=\{\tilde{Y}(t), t \in \R_+\}$ is said to be operator self-similar (o.s.s.) if
it is  continuous in law at each $t\geq 0$, and there exists $D\in End(\R^d)$ and nonrandom vectors $\{u(t), t\in\R_+\}$ in $\R^d$ such that
\beqnn
\big\{\tilde{Y}(ct)\big \}\stackrel{\mathscr{D}}{=}\big\{c^D \tilde{Y}(t)+u(c)\big\}\;\textrm{for all }\;c>0,
\eeqnn
where $\stackrel{\mathscr{D}}{=}$ denotes the equality of all finite-dimensional distributions, and $$
c^D=\exp \big((\log c)D\big)=\sum_{k=0}^{\infty} \frac{1}{k!} (\log c)^k D^k.
$$
The linear operator $D$  is called an {\it exponent} of the o.s.s. process $\tilde{Y}$.  For more information on this kind of processes, refer to  Cohen et al.\cite{CMR2010}, Hudson and Mason \cite{HM1982}, Laha and Rohatgi \cite{LR1982}, Marinucci and Robinson \cite{MR2000}, Meerschaert and Scheffler \cite[Chap.11]{MS2001}, and Sato \cite{ST1982}.

 Corresponding to the fractional Brownian motion (FBM) in one-dimensional case ($d=1$), there exists an operator
fractional Brownian motion (OFBM) in multidimensional case ($d\geq 2$). OFBMs are mean-zero, o.s.s., Gaussian processes with stationary increments.  They are of interest in several areas  for similar reasons  to  those in the univariate case.  For example, see
Chung \cite{Chung2002}, Davidson and de Jong \cite{Davidson2000}, Didier and Pipiras \cite{DP2011,DP2011a} and the references therein.

The asymptotical distribution of non-linear functionals of Gaussian vectors has been extensively studied. For example,  Arcones \cite{A1994}  considered limit theorems for functions of a stationary  Gaussian  sequence of vectors, and showed that the limit law can be either Gaussian or the law of a multiple Ito-Wiener integral, depending on the rate of decay of the coefficients. S\'{a}nchez \cite{S1993,S1995} studied limit theorems for non-linear functions of Gaussian vectors.  Inspired by these works, we  are also interested in this topic, which is the direct motivation of our work.

On the other hand, we should point out that Taqqu \cite{T1975} showed that the FBM can be approximated in law by a sequence of non-linear functions of Gaussian random variables. Noting that  OFBMs are the natural multivariate generalizations of FBMs, we are interested in whether  the OFBM can also be approximated in law by a sequence of non-linear functionals of  Gaussian vectors. Hence, we study limit theorems for   functionals of Gaussian vectors in this paper.

At the end of this section, we point out that all processes  considered here
are assumed to be proper.  We say that a process $\{X(t),t > 0\} $ is proper if for each $t > 0$ the distribution of $X(t)$ is full; that is, the
distribution is not contained in a proper hyperplane.

The rest of this paper is organized as follows. Section 2 is devoted to discussing weak convergence of  stationary $\R^d$-valued processes. In Section 3, we discuss  weak limit theorems for  functionals of Gaussian vectors. In Section 4,  we present an application of our results, and  show that a kind of OFBMs can be approximated in law by a sequence of functionals of Gaussian vectors.

\noindent
\\[4mm]

\noindent{\bbb 2\quad Sufficient conditions for weak convergence}\\[0.1cm]
 Let  $\big\{Z_N(t),t\in [0,\;1]\big\}_{N\in\N}$ be a sequence of $\R^d$-valued processes. In this section, we discuss the weak convergence of this sequence.  Before we state the main result of this section, we recall some basic facts. Throughout this paper,  let $B^*$ be the adjoint operator of $B\in End(\R^d)$, and $B^{-1}$ be the inverse of $B$.  We  use $\|x\|_2$ to denote the usual Euclidean norm of $x=(x^{(1)},\cdots,x^{(d)})^T \in \R^d$,  where $y^T$ denotes the transpose  of  $y\in\R^d$. Moreover, let
 $\left\|A\right\|=\max_{\|x\|_2=1}\|Ax\|_2$ denote the operator norm of $A\in End(\R^d)$. It is well-known that for any $A,B\in End({\R^d})$,
 \beqnn
 \left\|AB\right\|\leq \left\|A\right\|\cdot\left\|B\right\|,
 \eeqnn
and for every $A=(A_{ij})_{d\times d}\in End(\R^d)$,
\beqlb\label{s2-2}
\max_{1\leq i,j\leq d}|A_{ij}|\leq \left\|A\right\|\leq d^{\frac{3}{2}} \max_{1\leq i,j\leq d}|A_{ij}|.
\eeqlb
Furthermore, let
\beqnn
\lambda_A=\min\{\textrm{Re}\lambda: \lambda\in\sigma (A)\}\; \hbox{and}\;\Lambda_A=\max\{\textrm{Re}\lambda: \lambda\in\sigma(A)\},
\eeqnn
where $\sigma (A)$  is  the collection of all eigenvalues of $A$.

In order to state our results, we need to study the relationship between two linear operators  on $\R^d$.  For  any $n\in\N$, let $A(n)=\big(A_{ij}(n)\big)_{d\times d}\in End(\R^d)$ and $B(n)=\big(B_{ij}(n)\big)_{d\times d}\in End(\R^d)$.  We introduce the following  asymptotic notation.  We first introduce the small oh notation and the asymptotic equivalence.
\begin{defn}
Suppose that for any $i,j\in\{1,\cdots,d\}$, one of the following  cases holds.
\begin{itemize}
\item[(i)]There exists $N_0\in\N$ such that for all $n\geq N_0$,
\beqlb\label{a-17}
B_{ij}(n)\neq 0\;\textrm{and}\;\lim_{n\to\infty}A_{ij}(n)/B_{ij}(n)= a,
\eeqlb
where $a\in\R$.
\item[(ii)] There exists $N_1\in\N$ such that for all $n\geq N_1$,
\beqnn
A_{ij}(n)=0\;\textrm{and}\; B_{ij}(n)=0.
\eeqnn
\end{itemize}
If  $a=1$ in \eqref{a-17}, then we say that $A(n)$ is asymptotically equivalent to $B(n)$, as $n\to\infty$.
We denote this by $A(n)\sim B(n)$ as $n\to\infty$. If $a=0$ in \eqref{a-17}, then we say that $A(n)$ is of smaller order than $B(n)$, as $n\to\infty$.  We denote this by $A(n)= o(B(n))$ as $n\to\infty$.
\end{defn}
We have the following property.
\begin{lem}\label{rem4} If
$A(n)= o(B(n))\;\textrm{as}\;n\to\infty,
$
then
there exists an integer $N_0\in\N$  and a constant $K>0$ such that for all $n\geq N_0$,
\beqnn
\big\|A(n)\big\|\leq K\big\|B(n)\big\|.
\eeqnn

\end{lem}
The lemma \ref{rem4} can be easily proved.  Here we omit the proof.

 Next we introduce the big oh notation.

\begin{defn}  We write $A(n)=O\big(B(n)\big)$ as $n\to\infty$, if  there exists a constant $K> 0$ and an integer $N_0\in\N$ such that for all $n\geq N_0$,
\beqnn
| A_{ij}(n)|\leq K |B_{ij}(n)|\;\textrm{for all}\; i,j=1,\cdots,d.
\eeqnn

\end{defn}
We have the following property.
\begin{lem}\label{lem4}
If $A(n)=O\big(B(n)\big)$ as $n\to\infty$, then
 there exists an integer $N_0\in \N$ and a constant $K>0$ such that for all $n\geq N_0$,
\beqnn
\|A(n)\|\leq K \|B(n)\|.
\eeqnn

\end{lem}
It is easy to verify  that Lemma \ref{lem4} holds. Here we omit the proof.

\begin{defn}
Let $A=\big(A_{ij}\big)_{d\times d}\in End(\R^d) $ and $B=\big(B_{ij}\big)_{d\times d}\in End(\R^d)$. If
\beqnn
|A_{ij}|\leq |B_{ij}|\;\textrm{for all}\;i,j=1,\cdots,d,
\eeqnn  then we say $A\leq B$.
\end{defn}

We next introduce some technical lemmas which play an important role in our work.  The following lemma can be found in Mason and Xiao \cite{MX1999}.
\begin{lem}\label{lem1}
Let $D\in End(\R^d) $. If $\lambda_D>0$
and $r>0$, then for any $\delta>0$, there exist positive constants $K_1$ and $K_2$ such that
\beqnn
\left\|r^D\right\| \leq \begin{cases}K_1 r^{\lambda_D-\delta},  &\textrm{for all}\;  r\leq 1,
\\K_2 r^{\Lambda_D+\delta}, &\textrm{ for all}\; r\geq 1.
\end{cases}
\eeqnn
\end{lem}

In order to prove weak convergence, we need the following tightness criterion in the space $\mathcal{D}^d \big([0,\;1]\big)=\mathcal{D}^d\Big([0,\;1],\;\R^d\Big)$, which can be found in Dai \cite{Dai2011}.
\begin{lem}\label{lem2}
Let $\{Z_n(t), t\in[0,\;1]\}_{n\in\N}$ be a sequence of stochastic processes in $\mathcal{D}^d\big([0,\;1]\big)$ satisfying:
\begin{itemize}
\item[(i)] For every $n\in\N$, $Z_n(0)=0$ a.s.
\item[(ii)] There exist constants $K>0$, $\beta>0$,  $\alpha>1$ and an integer $N_0\in\N$ such that
\beqnn
\E\bigg[\Big\|Z_n(t)-Z_n(s)\Big\|_2^\beta\bigg] \leq K(t-s)^\alpha, n\geq N_0\;\textrm{and}\; 0\leq s\leq t\leq 1.
\eeqnn
\end{itemize}
Then $\{Z_n(t)\}$ is tight in $\mathcal{D}^d\big([0,\;1]\big)$.
\end{lem}

\medskip

 Let  $\{Y_i\}_{i\in\N}$  be  a stationary mean-zero sequence of random vectors with  $\E\big[\|Y_i\|_2^2\big]<\infty$. For any $N\in\N$, define
$$S_{\left\lfloor Nt\right\rfloor}=\sum_{i=1}^{\left\lfloor Nt\right\rfloor} Y_i,$$
where $\left\lfloor x\right\rfloor$ denotes the greatest
integer less than or equal to $x$. For convenience, let
\beqlb\label{FMC-2}S_N=\sum_{i=1}^{N} Y_i.\eeqlb
Furthermore, we assume that {\em empty sums are equal to $(0,\cdots,0)^T\in\R^d.$}

In the rest of this paper, most of  estimates  contain unspecified constants.
An unspecified positive and finite constant will be denoted by $\tilde{K}$,
which may not be the same in each occurrence. Sometimes we shall
emphasize the dependence of these constants upon parameters.   Moreover,  let $\Gamma $ denote a $d\times d$ symmetric and positive semi-definite matrix in the rest of this paper.

The main result of this section is the following.

\begin{lem}\label{thm1}
Suppose that a sequence $\{Z_N(t),t\in[0,\;1]\}_{N\in\N}$ of random functions  in  $\mathcal{D}^d([0,\, 1])$ satisfies:
\begin{itemize}
\item[(i)]\beqnn
Z_N(t)=N^{-D}B^{-1}S_{\left\lfloor Nt\right\rfloor},
\eeqnn
where  $D\in End(\R^d)$ with $\frac{1}{2}<\lambda_D,\Lambda_D<1$, and $B\in Aut(\R^d)$.
\item[(ii)]
\beqlb\label{thm1-a1}
\E\Big[S_N S^T_N\Big]=  BN^D \Gamma(N) N^{D^*}B^*,
\eeqlb
where $S_N$ is given by \eqref{FMC-2}, and $\Gamma(N)\in End(\R^d)$ with $\Gamma(N)\sim \Gamma$ as $N\to\infty$.

\item[(iii)]The finite-dimensional distributions of $\{Z_N(t)\}$ converge as $N\to\infty$.
\end{itemize}
Then the sequence $\{Z_N(t), t\in[0,\;1]\}$ converges weakly, as $N\to\infty$ in $\mathcal{D}^d([0,\;1])$, to an operator self-similar  process $X=\{X(t), t\in[0,\;1]\}$ with stationary increments, whose finite-dimensional distributions are the limits of those of $\{Z_N(t),t\in[0,\;1]\}$.
\end{lem}
{\it Proof of Lemma \ref{thm1}:} We choose $0\leq s\leq t\leq 1$. In order to prove Lemma \ref{thm1}, we first prove that $\{Z_N(t)\}$ is tight. In fact, we have
\beqlb\label{thm1-1}
\E\Big[\big\|Z_N(t)-Z_N(s)\big\|_2^{2}\Big]=\E\Big[\big\|Z_{\left\lfloor Nt\right\rfloor-\left\lfloor Ns\right\rfloor}\big\|_2^{2}\Big],
\eeqlb
since  $\{Y_i\}_{i\in\N}$ is stationary.

On the other hand, we  note that for any $x=(x^{(1)},\cdots,x^{(d)})^T\in\R^d$
\beqlb\label{thm1-9}
\|x\|_2^2=\sum_{k=1}^d(x^{(k)})^2.
\eeqlb
Hence, it follows from \eqref{thm1-1} and \eqref{thm1-9} that
\beqlb\label{thm1-10}
\E\Big[\big\|Z_{\left\lfloor Nt\right\rfloor-\left\lfloor Ns\right\rfloor}\big\|_2^{2}\Big]\leq \tilde{K}\bigg\|\E\Big[Z_{\left\lfloor Nt\right\rfloor-\left\lfloor Ns\right\rfloor}Z^T_{\left\lfloor Nt\right\rfloor-\left\lfloor Ns\right\rfloor}\Big]\bigg\|.
\eeqlb
We get from \eqref{thm1-a1} and \eqref{thm1-10}  that there exists $N_0\in\N$ such that for all $N\geq N_0$
\beqlb\label{thm1-11}
\E\Big[\big\|Z_{\left\lfloor Nt\right\rfloor-\left\lfloor Ns\right\rfloor}\big\|_2^{2}\Big]\leq \tilde{K}\Big\|[\frac{\left\lfloor Nt\right\rfloor-\left\lfloor Ns\right\rfloor}{N}]^D  \Big\| \times \Big\|  [\frac{\left\lfloor Nt\right\rfloor-\left\lfloor Ns\right\rfloor}{N}]^{D^*}\Big\|.
\eeqlb
Hence, it follows from  \eqref{thm1-11} and Lemma \ref{lem1} that for any $0<\delta<\lambda_D-\frac{1}{2}$
\beqnn
\E\Big[\big\|Z_N(t)-Z_N(s)\big\|_2^{2}\Big]\leq \tilde{K}[\frac{\left\lfloor Nt\right\rfloor-\left\lfloor Ns\right\rfloor}{N}]^{2(\lambda_D-\delta)},
\eeqnn
since $t,s\in [0,\;1]$.

On the other hand, due to de Haan \cite{dh1970}, we have
\beqlb\label{thm1-13}
\lim_{N\to\infty}[\frac{\left\lfloor Nt\right\rfloor-\left\lfloor Ns\right\rfloor}{N}]^{2(\lambda_D-\delta)}=(t-s)^{2(\lambda_D-\delta)}
\eeqlb
holds uniformly for $t\,,s\in[0,\;1]$.  Hence, it follows from  \eqref{thm1-11} and \eqref{thm1-13} that, for any $0<\delta<\lambda_D-\frac{1}{2}$, there exists a constant $N_0\in\N$ such that for all $N\geq N_0$
\beqlb\label{thm1-6}
\E\Big[\big\|Z_N(t)-Z_N(s)\big\|_2^{2}\Big]\leq \tilde{K} (t-s)^{2(\lambda_D-\delta)}.
\eeqlb
Finally, it follows from Lemma \ref{lem2} and \eqref{thm1-6} that $\{Z_N(t)\}$ is tight.

The tightness and convergence of the finite-dimensional distributions ((iii) of Lemma \ref{thm1}) ensure the weak convergence of $\{Z_N(t)\}$ to some limiting process $X=\{X(t)\}$. Since $\{Y_i\}$ is stationary, $\{X(t)\}$ must have  stationary increments.

Next, we show operator self-similarity.  It is obvious that $Z_N(0)=(0,\cdots,0)^T$. Hence, $X(0)=(0,\cdots,0)^T$.  Noting that $\{X(t)\}$ has stationary increments, we can easily get that $\{X(t)\}$ is continuous in law.
On the other hand,  for every $s>0$, let
\beqnn
\tilde{Z}(st)=\begin{cases}0,  &\textrm{if}\;  s\in(0,\;1),
\\{\left\lfloor s\right\rfloor}^{-D}B^{-1}S_{\left\lfloor s t\right\rfloor}, &\textrm{ if}\; s\geq 1.
\end{cases}
\eeqnn
It follows from (iii) of Lemma \ref{thm1} that the finite-dimensional distributions of $\{\tilde{Z}(st)\}$ converge to those of $\{X(t)\}$, as $s\to\infty$. From Theorem 5 in Hudson and Mason \cite{HM1982}, we get that $\{X(t)\}$ is operator self-similar. \qed
\begin{rem}\label{rem2}
From Lemma \ref{rem4} and the proof of Lemma \ref{thm1}, we can get that the condition (ii) in Lemma \ref{thm1} can be replaced by the following condition ($\Pi$).

\noindent{($\Pi$):}
\beqnn
\E[S_NS_N^T]=S_1(N)+S_2(N),
\eeqnn
where
\beqnn
S_1(N)=BN^D\Gamma(N)N^{D^*}B^*
\eeqnn
and
\beqnn
S_2(N)=BN^DA(N)N^{D^*}B^*
\eeqnn
with $A(N)=o(A)$ as $N\to\infty$ for some $A\in End(\R^d)$.
\end{rem}
\begin{rem} The matrix
$\Gamma$ is the covariance matrix of the limiting random vector $X(1)$.
\end{rem}

\noindent\\[4mm]

\noindent{\bbb 3\quad Limit theorems for non-Linear functionals }\\[0.1cm]
The main aim of this section is to  discuss  limit theorems for non-linear functionals  of Gaussian random vectors.  We will focus on a stationary Gaussian  sequence of $\R^d$-valued  random vectors $X_i=(X_i^{(1)},\cdots,X_i^{(d)})^T$ with
\beqlb\label{a-9} \E[X_i]=(0,\cdots,0)^T
\eeqlb and
\beqlb\label{a-5}
\E[X_i^{(p)}X_i^{(q)}]=\left\{
\begin{array}{ccc}
1,&\;\textrm{if}\;p=q,
\\0,&\;\textrm{others}.
\end{array}
\right.
\eeqlb
Let $\gamma(k)=\gamma(i,i+k)=\E\big[X_iX_{i+k}^T\big]=\big(\gamma_{pq}(i,i+k)\big)_{d\times d}$
be the  covariance matrix.
We are interested in  what conditions can be imposed on a function $G$
and on the sequence of covariance matrices  $\gamma(k)$ such that $\sum_{i=1}^{\left\lfloor Nt\right\rfloor}G(X_i)$
converges weakly to a process, as $N\to\infty$.

In order to answer the preceding question, we first introduce the following notation.
Let $$H_l(x)=(-1)^l e^{\frac{x^2}{2}}\frac{d^l}{dx^l}e^{-\frac{x^2}{2}},\;l\in\N\cup\{0\}$$
be the Hermite polynomials, and $X=(X^{(1)},\cdots,X^{(d)})^T$ be the standard $d$-dimensional Gaussian vector. For some fixed $L_i=\big(l^{(1)}_{i},\cdots,l^{(d)}_{i}\big)^T$, where  $i\in\{1,\cdots,d\}$ and $l_k^{(j)}\in \N\cup\{0\}$,  we
 define the following   $d$-dimensional random vector $e_{L_i}(X)$ by
\beqnn
e_{L_i}(X)=\bigg(e^{(1)}_{L_i}(X),\cdots,e^{(d)}_{L_i}(X)\bigg)^T,
\eeqnn
where the $j$th entry $e^{(j)}_{L_i}(X)$, $j=1,\cdots,d$, is given by
\beqnn
e^{(j)}_{L_i}(X)=\left\{
\begin{array}{ccc}
H_{l^{(1)}_1}(X^{(1)})\cdots H_{l_1^{(d)}}(X^{(d)}),&\;\textrm{if}\;i=j,
\\0,&\;\textrm{others},
\end{array}
\right.
\eeqnn
Furthermore, let $\mathscr{G}=\{G(x), x\in\R^d\}$ be the set of $\R^d$-valued measurable functions satisfying:
\begin{itemize}
\item[(i)] $\E[\|G(X)\|^2_2]<\infty$,
\item[(ii)] $\E[G(X)]=(0,\cdots,0)^T$.
\end{itemize}

Inspired by  Arcones \cite{A1994}, S\'{a}nchez \cite{S1993} and Taqqu \cite{T1975}, we define the following Hermite rank of a function $G\in\mathscr{G}$.
\begin{defn}\label{defn-2}Let $X$ be the standard $d$-dimensional  Gaussian vector and $G\in\mathscr{G}$. We define the {\it Hermite rank} of $G$ by
\beqnn
\textrm{{\em Rank}}\;(G)=\inf_{i\in\{1,\cdots,d\}}\Big\{\tau: \sum_{j=1}^dl_i^{(j)}=\tau\;\textrm{with}\; \E\Big[G^T(X)e_{L_i}(X)\Big]\neq 0\Big\}.\qquad
\eeqnn
\end{defn}
Moreover, $$\mathbb{G}_m=\big\{G: G\in\mathscr{G}\;\textrm{and} \;\textrm{{\it Rank}}\;(G)=m.\big\}.$$
\begin{rem}
From the definition \ref{defn-2}, we get that the rank of a function $G$ is unique. However, the corresponding index $L$ may be not unique.
\end{rem}
\begin{rem}
  The case that $Rank(G)=0$ is trivial, since $H_0(x)=1$. We will not discuss this trivial case. We assume that $Rank(G)\geq 1$ in the rest of this paper.
\end{rem}

\noindent\\[4mm]

\noindent{\bbb 3.1 \quad Conditions for weak convergence }\\[0.1cm]
In order to answer the problem in the previous part, we need  some additional conditions. Before we state them, we first introduce the following notation in the rest of this paper.
Let $\{X_i\}$ be a stationary mean-zero Gaussian sequence of $\R^d$-valued random vectors with \eqref{a-9} and \eqref{a-5}, and  $G\in\mathbb{G}_m$. Moreover, let $D\in End(\R^d)$ with $\frac{1}{2}<\lambda_D,\Lambda_D<1$, and $B\in Aut(\R^d)$.  For any $i,r \in\N$ and $n\in\{1,\cdots,d\}$, let
\beqnn
 I^{(n)}_{r,i}=\Big\{ L^{(i)}_n=(l_{(n,i)}^{(1)},\cdots, l_{(n,i)}^{(d)})^T:&&
 \sum_{j=1}^d l_{(n,i)}^{(j)}=r \nn
 \\&& \textrm{with}\;\E\big[G^T(X_{i})e_{L^{(i)}_{n}}(X_{i})\big]\neq 0\Big\}.
\eeqnn
Moreover, for any $i,j\in\N$ and $n_1,n_2\in\{1,\cdots,d\}$, we define  $I^{(n_1,n_2)}_{r,i,j}$ by
\beqnn
I^{(n_1,n_2)}_{r,i,j}=\Big\{(L^{(i)}_{n_1}, L^{(j)}_{n_2}):L^{(j)}_{n_2}\in I^{(n_2)}_{r,j}\;\textrm{and}\;L^{(i)}_{n_1}\in I^{(n_1)}_{r,i}\Big\}.
\eeqnn
At last, we use $E$ to denote the $d\times d$ matrix with all entries being $1$.

\begin{defn}\label{defn-1}  We say that $\{X_i\}$ satisfies  {\bf Condition} $\mathcal{H}(G,B,D,m)$ if
\begin{itemize}
\item[(i)] for some $\tilde{\Gamma}(N)=O\big(\Gamma \big)$ as $N\to\infty$,
\beqnn
\sum_{i,j=1}^{N}\Big(\sum_{p=1}^d\sum_{q=1}^d|\gamma_{pq}(i,j)|\Big)^m E =BN^D\tilde{\Gamma}(N)N^{D^*}B^*,
\eeqnn
where  $\gamma(i,j)=\big(\gamma_{pq}(i,j)\big)_{d\times d}$ is the covariance matrix given by
    \beqnn
    \gamma(i,j)=\gamma(|i-j|)=\E[X_iX_j^T];
    \eeqnn

\item[(ii)]  as $ |i-j|\to\infty$, \beqlb\label{condition-2}\Big\|\gamma\big(|i-j|)\Big\|\to 0;\eeqlb
\item[(iii)]
\beqlb\label{a-21}
\sum_{i,j=1}^{N}E(G,i,j,m)=BN^D\Gamma(N)N^{D^*}B^*,
\eeqlb
where \beqlb\label{a-16}
\Gamma(N)\sim \Gamma\;\textrm{ as}\; N\to\infty,\eeqlb and $E(G,i,j,m)=\big(E_{pq}(G,i,j,m)\big)_{d\times d}$ is given by
\beqnn
&&E_{pq}(G,i,j,m)=\nonumber
\\&&\qquad\sum_{(L^{(i)}_{p},L^{(j)}_{q})\in I^{(p,q)}_{m,i,j}}C_{L^{(i)}_{p} }C_{L^{(j)}_{q}}\bigg[\E\Big[\Pi_{n=1}^dH_{l^{(n)}_{(p,i)}}(X_i^{(n)})H_{l_{(q,j)}^{(n)}}(X_j^{(n)})\Big]\bigg]
\eeqnn
with
\beqnn
 C_{L^{(i)}_p}=\frac{\E\big[G^T(X_i)e_{L_p^{(i)}}(X_i)\big]}{\Pi_{k=1}^d l_{(p,i)}^{(k)}!}\;\textrm{for}\; L^{(i)}_p\in I^{(p)}_{m,i}.
\eeqnn

\end{itemize}
\end{defn}

Under {\bf Condition} $\mathcal{H}(G,B,D,m)$, we have the following result.
\begin{lem}\label{lem3}
 If $\{X_i\}$  satisfies  {\bf Condition} $\mathcal{H}(G,B,D,m)$, then
\beqlb\label{lem3-a1}
\E\bigg[\Big(\sum_{i=1}^NG(X_i)\Big)\Big(\sum_{i=1}^NG^T(X_i)\Big)\bigg]=S_1(N)+S_2(N),
\eeqlb
where
\beqnn
S_1(N)=BN^D\Gamma(N)N^{D^*}B^*
\eeqnn
and
\beqnn
S_2(N)=BN^Do(A)N^{D^*}B^*
\eeqnn
for some $A\in End(\R^d)$.
\end{lem}

\noindent{\it Proof of Lemma \ref{lem3}:} Since $X_i$ is the standard $d$-dimensional Gaussian vector, we can expand $G(X_i)$ as
\beqlb\label{lem3-2}
G(X_i)=\sum_{r\in\N}\Bigg\{\sum_{n=1}^d\sum_{L^{(i)}_n\in I^{(n)}_{r,i}}\big[C_{L^{(i)}_n} e_{L^{(i)}_n}(X_i)\big]\Bigg\},
\eeqlb
where $C_{L^{(i)}_n}=\frac{\E\big[G^T(X_i)e_{L_n^{(i)}}(X_i)\big]}{\Pi_{k=1}^d l_{(n,i)}^{(k)}!}$ if $L_n^{(i)}$ exists.  Moreover, if there exists some $n\in\{1,\cdots,d\}$ and $i\in\N$ such that $I_{r,i}^{(n)}=\emptyset$,  where $\emptyset$ denotes the null set,  then we assume that $\sum_{L^{(i)}_n\in I^{(n)}_{r,i}}\big[C_{L^{(i)}_n} e_{L^{(i)}_n}(X_i)\big]=(0,\cdots,0)^T$.

It follows from \eqref{lem3-2} and Sa\'nchez \cite{S1993} that for any $i,j\in\N$,
\beqlb\label{lem3-3}
&&\E\Big[G(X_i)G^T(X_j)\Big]\nonumber
\\&&=\E\bigg\{\sum_{r\in\N}\sum_{n_1,n_2=1}^d\sum_{(L_{n_1}^{(i)},L^{(j)}_{n_2})\in I^{(n_1,n_2)}_{r,i,j}} C_{L^{(i)}_{n_1}} C_{L^{(j)}_{n_2}} \Big[e_{L^{(i)}_{n_1}}(X_i)e^T_{L^{(j)}_{n_2}}(X_j)\Big]\bigg\}.
\eeqlb

Since $G\in\mathbb{G}_m$, we can rewrite the equation (\ref{lem3-3}) as follows.
\beqnn
\E\Big[G(X_i)G^T(X_j)\Big]=\E[\tilde{Q}(i,j)]+\E[\hat{Q}(i,j)],
\eeqnn
where
\beqlb\label{lem3-6}
\tilde{Q}(i,j)=\sum_{n_1,n_2=1}^d\tilde{Q}_{n_1n_2}(i,j)
 \eeqlb
 with
 $$\tilde{Q}_{n_1n_2}(i,j)=\sum_{(L^{(i)}_{n_1},L^{(j)}_{n_2})\in I^{(n_1,n_2)}_{m,i,j}}\Big[C_{L^{(i)}_{n_1}} C_{L^{(j)}_{n_2}} e_{L^{(i)}_{n_1}}(X_i)e^T_{L^{(j)}_{n_2}}(X_j)\Big],$$
and \beqlb\label{lem3-18-c}
\hat{Q}(i,j)=\sum_{n_1,n_2=1}^d\hat{Q}_{n_1n_2}(i,j)
\eeqlb with
$$\hat{Q}_{n_1n_2}(i,j)=\sum_{r\geq m+1\;\textrm{and}\;r\in\N}\sum_{(L^{(i)}_{n_1},L^{(j)}_{n_2})\in I^{(n_1,n_2)}_{r,i,j}}\Big[C_{L^{(i)}_{n_1}} C_{L^{(j)}_{n_2}} e_{L^{(i)}_{n_1}}(X_i)e^T_{L^{(j)}_{n_2}}(X_j)\Big].$$
Hence
\beqlb\label{lem3-5-b}
&&\E\bigg[\big(\sum_{i=1}^NG(X_i)\big)\big(\sum_{i=1}^NG^T(X_j)\big)\bigg]= \nonumber
\\&&\qquad\qquad\qquad\qquad\qquad \E\Big[\sum_{i=1}^N\sum_{j=1}^N \tilde{Q}(i,j)\Big]+\E\Big[\sum_{i=1}^N\sum_{j=1}^N \hat{Q}(i,j)\Big].
\eeqlb

In order to show \eqref{lem3-a1}, we first show that

\beqlb\label{lem3-5}
\E\bigg[\sum_{i=1}^N\sum_{j=1}^N \tilde{Q}(i,j)\bigg]=BN^D \Gamma(N) N^{D^*}B^*.
\eeqlb

By \eqref{lem3-6}, in order to show \eqref{lem3-5}, we need to focus on
\beqnn
\E\bigg[\sum_{i=1}^N\sum_{j=1}^N \tilde{Q}_{n_1n_2}(i,j)\bigg]\;\textrm{for all}\; n_1,n_2=1,\cdots,d.
\eeqnn
Here, we only look at the case that $n_1=n_2=1$. The other cases can be done in the same way.

\beqlb\label{lem3-8}
\E\Big[\tilde{Q}_{11}(i,j)\Big]=\E\bigg[\Big(M_{pq}(i,j)\Big)_{d\times d}\bigg]=\E[\mathbb{M}(i,j)],
\eeqlb
where
$$
M_{pq}(i,j)=\left\{
\begin{array}{ccc}
\sum_{(L^{(i)}_1,L^{(j)}_1)\in I^{(1,1)}_{m,i,j}}C_{L^{(i)}_{1} }C_{L^{(j)}_1}\Pi_{n=1}^dH_{l^{(n)}_{(1,i)}}(X_i^{(n)})H_{l_{(1,j)}^{(n)}}(X_j^{(n)}),&\;\textrm{if}\;p=q=1,
\\0,&\;\textrm{others}.
\end{array}
\right.
$$

Hence, we can get that
\beqlb\label{lem3-33}
\E\bigg[\sum_{i=1}^N\sum_{j=1}^N \tilde{Q}_{11}(i,j)\bigg]&&=
\sum_{i=1}^N\sum_{j=1}^N\Bigg(\sum_{(L^{(i)}_1,L^{(j)}_1)\in I^{(1,1)}_{m,i,j}}C_{L^{(i)}_{1} }C_{L^{(j)}_1}\nonumber
\\&&\bigg[\E\Big[\Pi_{n=1}^dH_{l^{(n)}_{(1,i)}}(X_i^{(n)})H_{l_{(1,j)}^{(n)}}(X_j^{(n)})\Big]\bigg]\Bigg) A(1,1),
\eeqlb
where \beqlb\label{lem3-33-e}
A(1,1)=\bigg[\begin{array}{cccc}1,&0,&\cdots,&0
\\ \quad &\quad&\cdots&\quad
\\0,&0,&\cdots,& 0 \end{array}\bigg]_{d\times d}.
\eeqlb

By using the same method as the proof of \eqref{lem3-33}, we can get that for any $n_1,n_2\in\{1,\cdots,d\}$
\beqlb\label{lem3-33-b}
&&\E\bigg[\sum_{i=1}^N\sum_{j=1}^N \tilde{Q}_{n_1n_1}(i,j)\bigg]=
\sum_{i=1}^N\sum_{j=1}^N\Bigg(\sum_{(L^{(i)}_{n_1},L^{(j)}_{n_2})\in I^{(n_1,n_2)}_{m,i,j}}C_{L^{(i)}_{n_1} }C_{L^{(j)}_{n_2}}\nn
\\&&\qquad\qquad\qquad\bigg[\E\Big[\Pi_{n=1}^dH_{l^{(n)}_{(n_1,i)}}(X_i^{(n)})H_{l_{(n_2,j)}^{(n)}}(X_j^{(n)})\Big]\bigg]\Bigg)A(n_1,n_2),
\eeqlb
where $A(n_1,n_2)=\big(A_{pq}(n_1,n_2)\big)_{d\times d}$ is a $d\times d$ matrix with
$$
A_{pq}(n_1,n_2)=\left\{
\begin{array}{ccc}
1,&\;\textrm{if}\;p=n_1\;\textrm{and}\;q=n_2,
\\0,&\;\textrm{others}.
\end{array}
\right.
$$
It follows from  (iii) of {\bf Condition} $\mathcal{H}(G,B,D,m)$, \eqref{lem3-33} and \eqref{lem3-33-b} that
\beqlb\label{lem3-33-c}
\E\bigg[\sum_{i=1}^N\sum_{j=1}^N \tilde{Q}(i,j)\bigg]=BN^D\Gamma(N)N^{D^*}B^*.
\eeqlb

By \eqref{lem3-5-b} and \eqref{lem3-5}, in order to establish (\ref{lem3-a1}),  we only need to show that
\beqlb\label{a-7}
\E\Big[\sum_{i=1}^N\sum_{j=1}^N \hat{Q}(i,j)\Big]=BN^Do(A)N^{D^*}B^*
\eeqlb
for some $A\in End(\R^d)$. To prove \eqref{a-7},
 we first show that as $N\to\infty$,
\beqlb\label{lem3-18-a}
\Bigg\|\E\bigg[N^{-D}B^{-1}\sum_{i=1}^N\sum_{j=1}^N \hat{Q}(i,j)(B^{*})^{-1}N^{-D^*}\bigg]\Bigg\|\to 0.
\eeqlb

By \eqref{lem3-18-c}, in order to prove \eqref{lem3-18-a},
we need  to look at the components $\hat{Q}_{n_1n_2}(i,j)$, $n_1,n_2=1,\cdots,d$.  On the other hand, according to (ii) of {\bf Condition} $\mathcal{H}(G,B,D,m)$, for arbitrarily small $\epsilon>0$, there exists $\tilde{Q}_0\in \N$ such that for all $|i-j|\geq \tilde{Q}_0$
\beqlb\label{lem3-32}
\sum_{p=1}^d\sum_{q=1}^d |\gamma_{pq}(|i-j|)|\leq \epsilon <1.
\eeqlb

In order to simplify the notation, let us define that for some integer $Q$ with $Q\geq\tilde{Q}_0$,
\beqnn
B(N,Q)&&=\{(i,j): |i-j|\leq Q, 0\leq i,j\leq N\}
\eeqnn
and
\beqnn
\tilde{B}(N,Q)&&=\{(i,j): |i-j|>Q, 0\leq i,j\leq N\}.
\eeqnn
Hence,  we have that for any $n_1,n_2\in\{1,\cdots,d\},$
\beqlb\label{lem3-15}
\sum_{i=1}^N\sum_{j=1}^N\Big[\hat{Q}_{n_1n_2}(i,j)\Big]
=\sum_{r\geq m+1\;\textrm{and}\;r\in\N}\tilde{V}_{n_1n_2}(Q)+\sum_{r\geq m+1\;\textrm{and}\;r\in\N}\hat{V}_{n_1n_2}(Q),
\eeqlb
where
$$
\tilde{V}_{n_1n_2}(Q)= \sum_{(i,j)\in B(N,Q)}\sum_{(L^{(i)}_{n_1},L^{(j)}_{n_2})\in I^{(n_1,n_2)}_{r,i,j}}\Big[C_{L^{(i)}_{n_1}} C_{L^{(j)}_{n_2}} e_{L^{(i)}_{n_1}}(X_i)e^T_{L^{(j)}_{n_2}}(X_j)\Big],
$$
and
$$
\hat{V}_{n_1n_2}(Q)=\sum_{(i,j)\in \tilde{B}(N,Q)}
\sum_{(L^{(i)}_{n_1},L^{(j)}_{n_2})\in I^{({n_1},{n_2})}_{r,i,j}}\Big[C_{L^{(i)}_{n_1}} C_{L^{(j)}_{n_2}} e_{L^{(i)}_{n_1}}(X_i)e^T_{L^{(j)}_n}(X_j)\Big].
$$
Therefore, we have
\beqnn
&&\Bigg\|\E\bigg[N^{-D}B^{-1}\sum_{i=1}^N\sum_{j=1}^N \hat{Q}(i,j)(B^*)^{-1}N^{-D^*}\bigg]\Bigg\|\nonumber
\\&&\quad\leq \sum_{n_1,n_2=1}^d \Bigg\|\E\bigg[N^{-D}B^{-1}\sum_{r\geq m+1\;\textrm{and}\;r\in\N} \Big[\tilde{V}_{n_1n_2}(Q)\Big](B^*)^{-1}N^{-D^*}\bigg]\Bigg\|\nonumber
\\&&\qquad+ \bigg\| \sum_{n_1,n_2=1}^d\E\Big[N^{-D}B^{-1}\sum_{r\geq m+1\;\textrm{and}\;r\in\N} \big[\hat{V}_{n_1n_2}(Q)\big](B^*)^{-1}N^{-D^*}\Big]\bigg\|.
\eeqnn

Next we  deal with \beqnn
\E\Big[N^{-D}B^{-1}\big[\sum_{r\geq m+1\;\textrm{and}\;r\in\N} \tilde{V}_{n_1n_2}(Q)\big](B^*)^{-1}N^{-D^*}\Big].
\eeqnn
We only focus on the case that $n_1=n_2=1$.   The other cases can be done in the same way.  The proof can be split into two steps.
We first assume that for any integer $r\geq m+1$ and $i,j\in\N$,
\beqlb\label{FMC-1}
I_{r,i,j}^{(1,1)}\neq \emptyset.
\eeqlb

Similar to \eqref{lem3-8}, we can get that
\beqlb\label{lem3-20}
\E\big[\sum_{r\geq m+1\;\textrm{and}\;r\in\N} \tilde{V}_{11}(Q)\big]=   \E\Big[\sum_{r\geq m+1\;\textrm{and}\;r\in\N} \sum_{(i,j)\in B(N,Q)} \mathcal{M}(i,j)\Big],
\eeqlb
where $\mathcal{M}(i,j)=\big(\mathcal{M}_{pq}(i,j)\big)_{d\times d}$ is a $d\times d$ matrix with
$$
\mathcal{M}_{pq}(i,j)=\left\{
\begin{array}{ccc}
\sum_{(L^{(i)}_1,L^{(j)}_1)\in I^{(1,1)}_{r,i,j}}C_{L^{(i)}_{1} }C_{L^{(j)}_1}\Pi_{n=1}^dH_{l^{(n)}_{(1,i)}}(X_i^{(n)})H_{l_{(1,j)}^{(n)}}(X_j^{(n)}),&\;\textrm{if}\;p=q=1,
\\0,&\;\textrm{others}.
\end{array}
\right.
$$
On the other hand, we have that
\beqlb\label{lem3-21}
&&\bigg|\sum_{r\geq m+1\;\textrm{and}\;r\in\N} \E[\mathcal{M}_{11}(i,j)]\bigg|\nonumber
\\&&=\bigg|\E\Big[\sum_{r\geq m+1\;\textrm{and}\;r\in\N} \sum_{(L_1^{(i)},L_1^{(j)})\in I^{(1,1)}_{r,i,j}}C_{L_1^{(i)}} C_{L_1^{(j)}} \prod_{n=1}^dH_{l_{(1,i)}^{(n)}}(X_i^{(n)})H_{l_{(1,j)}^{(n)}}(X_j^{(n)})\Big]\bigg|\nonumber
\\&&\leq \tilde{K} \Bigg|\E\bigg[\sum_{r\geq m+1\;\textrm{and}\;r\in\N} \sum_{L_1^{(i)}\in I^{(1)}_{r,i}}\Big[C_{L^{(i)}_1}\prod_{n=1}^dH_{l_{(1,i)}^{(n)}}(X_i^{(n)})\Big]^2\bigg]\Bigg|\nonumber
\\&&\qquad+\tilde{K} \Bigg|\E\bigg[\sum_{r\geq m+1\;\textrm{and}\;r\in\N} \sum_{L^{(j)}_1\in I^{(1)}_{r,j}}\Big[C_{L^{(j)}_1}\prod_{n=1}^dH_{l_{(1,j)}^{(n)}}(X_j^{(n)})\Big]^2\bigg]\Bigg|\nonumber
\\&&\leq\tilde{K}\E[\|G(X)\|_2^2]<\infty,
\eeqlb
where $X$ is the standard $d$-dimensional Gaussian vector.

From \eqref{lem3-20} and \eqref{lem3-21}, we get that
\beqlb\label{lem3-24}
\Bigg\|\E\bigg[\sum_{r\geq m+1\;\textrm{and}\;r\in\N}\tilde{V}_{11}(Q)\bigg]\Bigg\|\leq \tilde{K}(Q) N
\eeqlb
for some constant $\tilde{K}(Q)$ depending on $Q$.

On the other hand, by Lemma \ref{lem1}, we get that for any  $0<\delta<\lambda_D-\frac{1}{2}$,
\beqlb\label{lem3-35}
&&\Bigg\|\E\bigg[N^{-D}B^{-1}\sum_{r\geq m+1\;\textrm{and}\;r\in\N}\tilde{V}_{11}(Q)(B^*)^{-1}N^{-D^*}\bigg]\Bigg\|\nonumber
\\&&\qquad\leq \tilde{K} N^{-2(\lambda_D-\delta)}\bigg\|\E\Big[\sum_{r\geq m+1\;\textrm{and}\;r\in\N}\tilde{V}_{11}(Q)\Big]\bigg\|.
\eeqlb
From \eqref{lem3-24} and \eqref{lem3-35}, we have
\beqlb\label{a-14}
\Bigg\|\E\bigg[N^{-D}B^{-1}\sum_{r\geq m+1\;\textrm{and}\;r\in\N}\tilde{V}_{11}(Q)(B^*)^{-1}N^{-D^*}\bigg]\Bigg\|\leq \tilde{K}(Q) N^{-2(\lambda_D-\delta)+1}.
\eeqlb
By using the same method as the proof of \eqref{a-14}, we can get that

\beqlb\label{a-13}
&&\sum_{n_1,n_2=1}^d\Bigg\|\E\bigg[N^{-D}B^{-1}\sum_{r\geq m+1\;\textrm{and}\;r\in\N}\tilde{V}_{n_1n_2}(Q)(B^*)^{-1}N^{-D^*}\bigg]\Bigg\|\nn
\\&&\qquad\qquad \leq \tilde{K}(Q) N^{-2(\lambda_D-\delta)+1}.
\eeqlb

Now we turn to \beqlb\label{a-11}
\sum_{n_1,n_2=1}^d\E\Big[N^{-D}B^{-1}\sum_{r\geq m+1\;\textrm{and}\;r\in\N}  \big[\hat{V}_{n_1n_2}(Q)\big](B^*)^{-1}N^{-D^*}\Big].
\eeqlb We first look at $ \hat{V}_{11}(Q)$. Similar  to \eqref{lem3-8}, we have
\beqlb\label{lem3-25}
\E\big[ \hat{V}_{11}(Q)\big]&&= \sum_{(i,j)\in \tilde{B}(N,Q)}\bigg(\E\big[\mathcal{M}(i,j)\big]\bigg).
\eeqlb
We also note that
\beqlb\label{lem3-27}
&&\Big|\E\big[\mathcal{M}_{11}(i,j)\big]\Big|\leq \sum_{(L_1^{(i)},L_1^{(j)})\in I^{(1,1)}_{r,i,j}} \Bigg| \E\bigg[C_{L^{(i)}_1 } C_{L^{(j)}_1}\prod_{n=1}^dH_{l_{(1,i)}^{(n)}}(X_i^{(n)})H_{l_{(1,j)}^{(n)}}(X_j^{(n)})\bigg]\Bigg|\nonumber
\\&&\quad=\sum_{(L_1^{(i)},L_1^{(j)})\in I^{(1,1)}_{r,i,j}}\bigg\{\Bigg| \E\bigg
[C_{L_1^{(i)}} C_{L^{(j)}_1}r!\prod_{n=1}^d\frac{H_{l_{(1,i)}^{(n)}}(X_i^{(n)})H_{l_{(1,j)}^{(n)}}(X_j^{(n)})}{l_{(1,i)}^{(n)}!l_{(1,j)}^{(n)}!}\bigg]\Bigg|\nonumber
\\&&\qquad  \frac{\prod_{n=1}^dl_{(1,i)}^{(n)}!l_{(1,j)}^{(n)}!}{r!}\bigg\}.
\eeqlb
By the Cauchy-Schwartz  inequality, we get that
\beqlb\label{lem3-26}
&&\sum_{(L_1^{(i)},L_1^{(j)})\in I^{(1,1)}_{r,i,j}}  \E\Bigg|\bigg[C_{L_1^{(i)}} C_{L_1^{(j)}}\prod_{n=1}^d H_{l_{(1,i)}^{(n)}}(X_i^{(n)})H_{l_{(1,j)}^{(n)}}(X_j^{(n)})\bigg]\Bigg| \nonumber
\\&& \leq \bigg\{\sum_{(L_1^{(i)},L_1^{(j)})\in I^{(1,1)}_{r,i,j}}\big(C_{L^{(i)}_1})^2 (C_{L^{(j)}_1})^2\Big(\frac{\prod_{n=1}^d l^{(n)}_{(1,i)}!l^{(n)}_{(1,j)}!}{r!}\Big)^{2}\bigg\}^{\frac{1}{2}}\nonumber
\\&&\quad\times \Bigg\{\sum_{(L_1^{(i)},L_1^{(j)})\in I^{(1,1)}_{r,i,j}}\Bigg(r!\E\big[\prod_{n=1}^d \frac{H_{l_{(1,i)}^{(n)}}(X_i^{(n)})H_{l_{(1,j)}^{(n)}}(X^{(n)}_j)}{l_{(1,i)}^{(n)}! l_{(1,j)}^{(n)}!}\big]\Bigg)^2\Bigg\}^{\frac{1}{2}}.
\eeqlb
On the other hand, we have that
\beqlb\label{lem3-28}
&&\Bigg\{\sum_{(L_1^{(i)},L_1^{(j)})\in I^{(1,1)}_{r,i,j}}\bigg(r!\Big|\E\big[\prod_{n=1}^d\frac{H_{l_{(1,i)}^{(n)}}(X_i^{(n)})H_{l_{(1,j)}^{(n)}}(X_j^{(n)})}{l_{(1,i)}^{(n)}!l_{(1,j)}^{(n)}!}\big]\Big|\bigg)^2\Bigg\}^{\frac{1}{2}}\nonumber
\\&&\quad\qquad\qquad \leq \tilde{K} \sum_{(L_1^{(i)},L_1^{(j)})\in I^{(1,1)}_{r,i,j}}\bigg|r!\E\Big[\prod_{n=1}^d \frac{H_{l_{(1,i)}^{(n)}}(X_i^{(n)})H_{l_{(1,j)}^{(n)}}(X_j^{(n)})}{l_{(1,i)}^{(n)}!l_{(1,j)}^{(n)}!}\Big]\bigg|.
\eeqlb
Moreover,  due to S\'{a}nchez \cite{S1993}, we obtain that
\beqlb\label{lem3-30}
&&\sum_{(L_1^{(i)},L_1^{(j)})\in I^{(1,1)}_{r,i,j}}\bigg|r!\E\Big[\prod_{n=1}^d \frac{H_{l_{(1,i)}^{(n)}}(X_i^{(n)})H_{l_{(1,j)}^{(n)}}(X_j^{(n)})}{l_{(1,i)}^{(n)}!l_{(1,j)}^{(n)}!}\Big]\bigg|\nn
\\&&\qquad\qquad\qquad\leq \Big(\sum_{p=1}^d\sum_{q=1}^d|\gamma_{pq}(i,j)|\Big)^r.
\eeqlb
Finally, we note that
\beqlb\label{lem3-29}
&&\bigg\{\sum_{(L_1^{(i)},L_1^{(j)})\in I^{(1,1)}_{r,i,j}}(C_{L^{(j)}_1})^2 (C_{L^{(i)}_1})^2\Big(\frac{\prod_{i=1}^nl_{(1,i)}^{(n)}!l_{(1,j)}^{(n)}!}{r!}\Big)^2\bigg\}^{\frac{1}{2}}\nonumber
\\&&\qquad\qquad\qquad\leq \tilde{K} \Big(\sum_{L^{(i)}_1\in I^{(1,1)}_{r,i}}(C_{L^{(i)}_1})^2 \Pi_{n=1}^dl_{(1,i)}^{(n)}!\Big),\qquad
\eeqlb
since $$\Pi_{n=1}^d l_{(1,j)}^{(n)}!\leq r!\;\textrm{and}\;\Pi_{n=1}^d l_{(1,i)}^{(n)}!\leq r!.$$  It follows from (\ref{lem3-25}) to (\ref{lem3-29}) that
\beqlb\label{lem3-31}
&&\sum_{r\geq m+1\;\textrm{and}\;r\in\N}\E \big[\hat{V}_{11}(Q)\big]\nonumber
\\&&\leq \tilde{K }
\sum_{r\geq m+1\;\textrm{and}\;r\in\N} \sum_{(i,j)\in \tilde{B}(N,Q)}
\Big(\sum_{p=1}^d\sum_{q=1}^d|\gamma_{pq}(i,j)|\Big)^r\nn
\\&&\qquad\qquad
\Big(\sum_{L^{(i)}_1\in I^{(1)}_{r,i}}(C_{L^{(i)}_1})^2\Pi_{n=1}^dl_{(1,i)}^{(n)}!\Big)A(1,1)\nonumber
\\&&\leq \tilde{K}\sum_{r\geq m+1\;\textrm{and}\;r\in\N}
 \sum_{(i,j)\in\tilde{B}(Q,N)}
 \Big(\sum_{p=1}^d\sum_{q=1}^d|\gamma_{pq}(|i-j|)|\Big)^r\nn
 \\&&\qquad\qquad
  \Big(\sum_{L^{(i)}_1\in I^{(1)}_{r,i}}(C_{L^{(i)}_1})^2\Pi_{n=1}^dl_{(1,i)}^{(n)}!\Big)A(1,1),\qquad\quad
\eeqlb
where $A(1,1)$ is given by \eqref{lem3-33-e}.

Let $$C_G(X)=\sum_{r=0}^\infty\sum_{L\in I_r}C_{L}^2 \Pi_{j=1}^dl^{(j)}!.$$
Then
\beqlb\label{a-12}
C_G(X)<\infty,
\eeqlb
since $\E\Big[\|G(X)\|_2^2\Big]<\infty$.

By \eqref{lem3-32}, \eqref{lem3-31} and \eqref{a-12}, we get that
\beqlb\label{lem3-33-a}
\sum_{r\geq m+1\;\textrm{and}\;r\in\N} \E \big[\hat{V}_{11}(Q)\big]\leq \tilde{K}  \epsilon  \sum_{(i,j)\in \tilde{B}(N,Q)}
\Big(\sum_{p=1}^d\sum_{q=1}^d|\gamma_{pq}(i,j)|\Big)^m  A(1,1).
\eeqlb
By  \eqref{lem3-33-a},
\beqlb\label{a-8}
\sum_{r\geq m+1\;\textrm{and}\;r\in\N} \E \big[\hat{V}_{11}(Q)\big]&&\leq \tilde{K}\epsilon  \sum_{i,j=1}^N
\Big(\sum_{p=1}^d\sum_{q=1}^d|\gamma_{pq}(i,j)|\Big)^m  A(1,1).
\eeqlb
By using the same method as the proof of \eqref{a-8}, we get that
\beqnn
\sum_{n_1,n_2=1}^d\sum_{r\geq m+1\;\textrm{and}\;r\in\N} \E \big[\hat{V}_{n_1n_2}(Q)\big]&&\leq \tilde{K}\epsilon  \sum_{i,j=1}^N
\Big(\sum_{p=1}^d\sum_{q=1}^d|\gamma_{pq}(i,j)|\Big)^m  E.
\eeqnn
From (i) of {\bf Condition} $\mathcal{H}(G,B,D,m)$, we have that
\beqnn
\sum_{n_1,n_2=1}^d\sum_{r\geq m+1\;\textrm{and}\;r\in\N} \E \big[\hat{V}_{n_1n_2}(Q)\big]\leq \tilde{K}\epsilon BN^D\tilde{\Gamma}(N) N^{D^*}B^*.
\eeqnn
Then, there exists $\tilde{\Gamma}\in End(\R^d)$ such that
\beqlb\label{a-2}
\sum_{n_1,n_2=1}^d\sum_{r\geq m+1\;\textrm{and}\;r\in\N} \E \big[N^{-D}B^{-1}\hat{V}_{n_1n_2}(Q)(B^*)^{-1}N^{-D^*}\big]= o\big(\tilde{\Gamma}\big).
\eeqlb

From \eqref{lem3-32}, \eqref{a-13} and \eqref{a-2}, we get that as $N\to\infty$
\beqlb\label{lem3-36}
&&\sum_{n_1,n_2=1}^d\Bigg\|\E\bigg[N^{-D}B^{-1}\sum_{r\geq m+1\;\textrm{and}\;r\in\N} \tilde{V}_{n_1n_2}(Q)(B^*)^{-1}N^{-D^*}\bigg]\Bigg\|\nonumber
\\&&+\Bigg\|\sum_{n_1,n_2=1}^d\E\bigg[N^{-D}B^{-1}\sum_{r\geq m+1\;\textrm{and}\;r\in\N} \hat{V}_{n_1n_2}(Q)(B^*)^{-1}N^{-D^*}\bigg]\Bigg\|\to 0.
\eeqlb

Next we assume that \eqref{FMC-1} does not hold. It follows  from the above proof that \eqref{lem3-36} still holds.

Combining \eqref{lem3-15} and \eqref{lem3-36}, we get \eqref{lem3-18-a}.

Next we prove \eqref{a-7}. We first  point out that for any $A(n)\in End(\R^d)$, if
\beqnn
\|A(n)\|\to 0\;\textrm{as}\;n\to\infty,
\eeqnn
then for all $i,j=1,\cdots,d$,
\beqnn
A_{ij}(n)\to 0\;\textrm{as}\; n\to\infty.
\eeqnn
From \eqref{lem3-18-a},  we get that there exists some  $A\in End(\R^d)$ such that
\beqnn
\E\bigg[N^{-D}B^{-1}\sum_{i=1}^N\sum_{j=1}^N \hat{Q}(i,j)(B^*)^{-1}N^{-D^*}\bigg]= o\big(A\big)\;\textrm{as}\; N\to\infty.
\eeqnn
Then we get that
\beqlb\label{a-6}
\E\bigg[\sum_{i=1}^N\sum_{j=1}^N \hat{Q}(i,j)\bigg]=BN^{D} o\big(A\big)N^{D^*}B^*.
\eeqlb
From \eqref{lem3-33-c} and \eqref{a-6}, we  get that the lemma holds.
\qed

\begin{rem}\label{rem1}
From the proof of Lemma \ref{lem3}, we easily get that as $N\to\infty$,
\beqlb\label{rem1-1}
&&\Bigg\|\sum_{n_1,n_2=1}^d\sum_{i,j=1}^N\E\bigg[\sum_{r\geq m+1\;\textrm{and}\;r\in\N}\sum_{(L^{(i)}_{n_1},L^{(j)}_{n_2})\in I^{(n_1,n_2)}_{r,i,j}}\nonumber
\\&&N^{-D}B^{-1}\Big[C_{L^{(i)}_{n_1}} C_{L^{(j)}_{n_2}} e_{L^{(i)}_{n_1}}(X_i)e^T_{L^{(j)}_{n_2}}(X_j)\Big](B^{*})^{-1}N^{-D^*}
\bigg]\Bigg\|\to 0.
\eeqlb
\end{rem}

\noindent\\[4mm]

\noindent{\bbb 3.2 \quad The reduction theorem }\\[0.1cm]
In this subsection, we assume that $G\in\mathbb{G}_m$, and $\{X_i\}$ satisfies {\bf Condition} $\mathcal{H}(G,B,D,m)$. We  study  weak limit theorems for the process
\beqlb\label{thm2-a4}
Z_N(t)=N^{-D}B^{-1}\sum_{i=1}^{\left\lfloor Nt \right\rfloor}G(X_i),\;t\in[0,1].
\eeqlb

For any $t\in[0,\;1]$, define
\beqlb\label{thm2-a1}
Z_{N,m}(t)&&=N^{-D}B^{-1}\bigg[\sum_{i=1}^{\left\lfloor Nt\right\rfloor}\Big[\sum_{n=1}^d\sum_{L^{(i)}_n\in I^{(n)}_{m,i}}C_{L^{(i)}_n} e_{L^{(i)}_n}(X_i)\Big]\bigg],
\eeqlb
and
\beqlb\label{thm2-a2}
\tilde{Z}_{N,m}(t)=N^{-D}B^{-1}\sum_{i=1}^{\left\lfloor Nt\right\rfloor}\sum_{r\geq m+1\;\textrm{and}\;r\in\N}\bigg[\sum_{n=1}^d\sum_{L^{(i)}_n\in I^{(n)}_{r,i}}C_{L^{(i)}_n} e_{L^{(i)}_n}(X_i)\bigg].
\eeqlb

Before we state our result, we need the following useful lemma.

\begin{lem}\label{thm2}
If the limit in distribution of $\big(Z_{N,m}(t_1),\cdots,Z_{N,m}(t_p)\big)$ exists $\Big($we denote it by $\big(Z_m(t_1),\cdots,Z_m(t_p)\big)\Big)$, then as $N\to\infty$
\beqlb\label{thm3-1}
\Big(Z_{N}(t_1),\cdots, Z_{N}(t_p)\Big)\stackrel{\mathscr{D}}{\Rightarrow}\bigg( Z_m(t_1),\cdots,Z_m(t_p)\bigg),\quad t_1,\cdots,t_p\in[0,\;1],
\eeqlb
where $\stackrel{\mathscr{D}}{\Rightarrow}$ denotes convergence in distribution.
\end{lem}
{\it Proof of Lemma \ref{thm2}:} In order to simplify the discussion,  we only prove the case that  $p=1$. The general case can be done in the same way. According to \eqref{lem3-2}, we have
\beqlb\label{a-1}
Z_N(t)=N^{-D}B^{-1}\sum_{i=1}^{\left\lfloor Nt\right\rfloor}\sum_{r\geq m\;\textrm{and}\;r\in\N}\sum_{n=1}^d\sum_{L^{(i)}_n\in I^{(n)}_{r,i}}C_{L^{(i)}_n}e_{L^{(i)}_n}(X_i).
\eeqlb
From  \eqref{thm2-a1}, \eqref{thm2-a2} and \eqref{a-1},  we have
\beqnn
Z_N(t)=Z_{N,m}(t)+\tilde{Z}_{N,m}(t).
\eeqnn
Hence, in order to prove (\ref{thm3-1}), it is sufficient to prove that $\{\tilde{Z}_{N,m}(t)\}$ converges to the $d$-dimensional zero vector in probability, that is,  as $N\to\infty$
\beqlb\label{thm2-3}
\P\Big\{\big\|\tilde{Z}_{N,m}(t)\big\|_2\geq \epsilon \Big\}\to 0.
\eeqlb

Note that for an $\R^d$-valued random variable $Y=(Y^{(1)},\cdots,Y^{(d)})^T$, $\E[\|Y\|_2^2]$ equals the sum of diagonal entries of the covariance matrix $\E[YY^T]$. It follows from  \eqref{s2-2} and \eqref{thm2-a2} that
\beqlb\label{lem3-39}
\E\Big[\Big\|\tilde{Z}_{N,m}(t)\Big\|_2^2\Big]&&\leq \tilde{K} \bigg\|\E\Big[\tilde{Z}_{N,m}(1)\tilde{Z}^T_{N,m}(1)\Big]\bigg\|,
\eeqlb
since $t\in [0,\;1]$.
By \eqref{lem3-39},
\beqlb\label{thm2-4}
&&\E\Big[\big\|\tilde{Z}_{N,m}(t)\big\|_2^2\Big]\nonumber
\\&&\qquad\leq \tilde{K}\Bigg\|\sum_{n_1,n_2=1}^d\E\bigg[ N^{-D}B^{-1}\sum_{i,j=1}^N\sum_{r\geq m+1\;\textrm{and}\;r\in\N}
\sum_{(L^{(i)}_{n_1},L^{(j)}_{n_2})\in I^{(n_1,n_2)}_{r,i,j}}\nonumber
\\&&\qquad\qquad C_{L^{(i)}_{n_1}} C_{L^{(j)}_{n_2}}  \big[e_{L^{(i)}_{n_1}}(X_i)e^T_{L^{(j)}_{n_2}}(X_j)\big](B^*)^{-1}N^{-D^*}\bigg]\Bigg\|.
\eeqlb
Therefore, we get from Remark \ref{rem1}   and the Chebyshev-Markov inequality \cite[Chap.1]{David2004} that (\ref{thm2-3}) holds. So the lemma holds. \qed

\begin{thm}\label{thm}
Let $G\in \mathbb{G}_m$ for some $m\in\N$, and $\{X_i\}$ satisfy {\bf Condition} $\mathcal{H}(G,B,D,m)$. Define $Z_N(\cdot)$ as in (\ref{thm2-a4}) and $Z_{N,m}(\cdot)$ as in (\ref{thm2-a1}). If the finite-dimensional distributions of $\{Z_{N,m}(\cdot)\}$ converge to those of some process $\{Z_m(\cdot)\}$, then $\{Z_N(\cdot)\}$ converges weakly to the process $\{Z_m(\cdot)\}$ in $\mathcal{D}^d([0,1])$.
\end{thm}

\noindent{\it Proof of Theorem \ref{thm}:} In order to prove the theorem, it suffices to prove that $\{Z_N(t)\}$  satisfies Lemma \ref{thm1}.  By Lemma \ref{lem3},  the condition ($\Pi$) in Remark \ref{rem2} holds. Finally,  Lemma \ref{thm2} implies that  the  condition (iii) in Lemma \ref{thm1} is satisfied.  From the above arguments, we get that the theorem holds. \qed

\noindent\\[4mm]

\noindent{\bbb 4 \quad Application }\\[0.1cm]
As an application of our results,   we show that, under some additional conditions, the limiting process of $\{Z_N(t),t\in[0,\;1]\}$ given by \eqref{thm2-a4} is, up to a multiplicative matrix from the left,  a time-reversible operator fractional Brownian motion.

We first recall an integral representation of OFBM. Let $D$ be a linear operator on $\R^d$ with $0<\Lambda_D,\lambda_D<1$.  Moreover, let  $X=\{X(t)\}$ be an OFBM with o.s.s. exponent $D$.  Then, from Didier and Pipiras \cite{DP2011}, we know that  $X$ admits the following integral representation
\beqnn
\{X(t)\}\stackrel{\mathscr{D}}{=}\Big\{\int_{\R}\frac{e^{itx}-1}{ix} \Big(x_+^{-(D-\frac{I}{2})}A+x_-^{-(D-\frac{I}{2})}\bar{A}\Big)W(dx)\Big\}
\eeqnn
for some linear  operator $A$ on $\mathbb{C}^d$. Here, $\bar{A}$ denotes the complex conjugate and
\beqnn
W(x):=W_1(x)+iW_2(x)
\eeqnn
denotes a complex-valued multivariate Brownian motion such that $W_1(-x)=W_1(x)$ and $W_2(-x)=-W_2(x)$,  $W_1(x)$ and $W_2(x)$ are independent, and the induced random measure $W(x)$ satisfies
$$
\E\Big[W(dx)W^*(dx)\Big]=dx,
$$
where $W^*$ is the adjoint operator of $W$.   Moreover, it follows from Dai \cite{Dai2011} that, up to a multiplicative constant, we can rewrite $\{X(t)\}$ as follows.
\beqlb\label{a-15}
\{X(t)\}\stackrel{\mathscr{D}}{=}\Big\{\int_0^\infty G_1(x,t)W_1(dx)+\int_0^\infty G_2(x,t)W_2(dx)\Big\},
\eeqlb
where
\beqnn
G_1(x,t)&&=\frac{\sin tx}{x} x^{-(D-\frac{I}{2})}A_1+\frac{\cos tx -1}{x} x^{-(D-\frac{I}{2})}A_2,
\\G_2(x,t)&&=\frac{\sin tx}{x} x^{-(D-\frac{I}{2})}A_2+\frac{1-\cos tx }{x} x^{-(D-\frac{I}{2})}A_1,
\eeqnn
and
\beqnn
A=A_1+iA_2.
\eeqnn

  In order to reach our aim in this section, we need the following technical lemma.
\begin{lem}\label{lem6}
Let $\big\{Z_i\big\}_{i\in\N}$ be  a stationary mean-zero Gaussian  sequence of $\R^d$-valued vectors. Let $$\tilde{\gamma}(i,j)=\E[Z_iZ_j^T].$$  Suppose that
\beqnn
\sum_{i=1}^N\sum_{j=1}^N \tilde{\gamma}(i,j)= \tilde{K} UN^D \Gamma_1(N) N^{D^*}U^*,
\eeqnn
where $U\in Aut(\R^d)$ and $\Gamma_1(N)\sim\Gamma_1$ as $N\to\infty$ with $\Gamma_1=\E[ X(1)X^T(1)]$.
Then,
\beqnn
Q_N(t)=N^{-D}U^{-1}\sum_{i=1}^{\left\lfloor Nt\right\rfloor}Z_i,
\eeqnn
converges weakly, as $N\to\infty$ in $\mathcal{D}^d([0,1])$, up to a multiplicative matrix from the left, to the time-reversible OFBM $X$  given by \eqref{a-15} with $A_2A_1^*=A_1A_2^*$.
\end{lem}

 By using the same method as the proof of Theorem 2.2 in Dai \cite{Dai2011}, we can easily prove  this lemma. Here we omit the proof.

Next, we state the main result of this section.
\begin{cor}\label{thm3}
Suppose that $\{X_i\}$ satisfies {\bf Condition} $\mathcal{H}(G,B,D,1)$ with $\Gamma=\Gamma_1$ in \eqref{a-16}.
Then,\beqnn
Z_N(t)=N^{-D}B^{-1}\sum_{i=1}^{\left\lfloor Nt\right\rfloor}G(X_i),\;t\in[0,\;1],
\eeqnn
converges weakly, as $N\to\infty$ in $\mathcal{D}^d([0,\;1])$, up to a multiplicative matrix from the left, to the time-reversible OFBM $X$ given by \eqref{a-15} with $A_2A_1^*=A_1A_2^*$.
\end{cor}
{\it Proof of Corollary \ref{thm3}:}  It follows from Lemma \ref{lem3} and Theorem \ref{thm} that, in order to prove Corollary \ref{thm3}, it suffices to show that $Z_{N,1}(t)$ given by \eqref{thm2-a1} converges weakly, as $N\to\infty$ in $\mathcal{D}^d([0,\;1])$, up to a multiplicative matrix from the left, to the time-reversible OFBM $X$.   In fact, since $Z_{N,1}(t)$  is proper and $H_1(x)=x$, we can get that there exists $C\in Aut(\R^d)$ such that
\beqnn
Z_{N,1}(t)=N^{-D}B^{-1}\sum_{i=1}^{\left\lfloor Nt\right\rfloor}CX_i.
\eeqnn
 Since $\{X_i\}$ is stationary and Gaussian, so is $\{CX_i\}$.

For convenience, let $\tilde{Z}_i=CX_i$.  Next, we check that $\{\tilde{Z}_i\}$ satisfies Lemma \ref{lem6}.
In fact, it follows from \eqref{a-21} that
\beqnn
\sum_{i,j=1}^N \E[\tilde{Z}_i\tilde{Z}^T_j]=BN^D\Gamma(N)N^{D^*}B^*.
\eeqnn
Hence, it follows from Lemma \ref{lem6}  that $Z_{N,1}(t)$ converges weakly, up to a multiplicative matrix from the left, to the time-reversible OFBM $X$. Finally, we get that the corollary holds.\qed
\medskip
\begin{rem} In Corollary \ref{thm3}, {\bf Condition} $\mathcal{H}(G,B,D,1)$ implies that $\frac{1}{2}<\lambda_D,\Lambda_D<1$. For OFBMs, the condition $\frac{1}{2}<\lambda_D,\Lambda_D<1$  in the univariate case is known as the long range dependence (LRD).  In the multivariate case, the condition has the potential to generate a divergence of the spectrum at zero. See Didier and Pipiras \cite{DP2011}.  Hence, we may define the operator LRD in the sense of  $\frac{1}{2}<\lambda_D,\Lambda_D<1$.  There is only a little work related to this topic. See, for example,  Didier and Pipiras \cite{DP2011}.  However, considering the importance of LRD in applications, it is worth  spending much more time on the LRD in the multivariate context.
\end{rem}
\noindent\\[4mm]
\noindent\bf{\footnotesize Acknowledgements}\quad\rm
{\footnotesize We  would
like to thank the reviewers for their helpful comments which greatly improve this work.}\\[4mm]

\noindent{\bbb{References}}


\begin{enumerate}
{\footnotesize
\bibitem{A1994}
Arcones  M A.
\newblock{Limit theorems for nonlinear funcitonals of a stationary Gaussian sequence of vectors.}
\newblock{The Annals of Probability,} 1994, 22: 2242--2272.\\[-6.5mm]
\bibitem{Chung2002}
Chung  C F.
\newblock{Sample means, sample autocovariances, and linear regression of stationary multivariate long memory processes.}
\newblock{Econometric Theory, } 2002, 18: 51--78.\\[-6.5mm]
\bibitem{CMR2010}
Cohen S, Meerschaert M M, Rosi\'nski J.
\newblock{Modeling and simulation with operator scaling.}
\newblock{Stoch Process Appl, }2010, 120: 2390--2411.\\[-6.5mm]
\bibitem{D1970}
Davydov  Y.
\newblock{The invariance principle for
stationary processes.}
\newblock{ Teor Verojatnost i Primenen,} 1970, 15: 498--509.\\[-6.5mm]
\bibitem{Davidson2000}
Davidson  J, De Jong R M.
\newblock{The functional central limit theorem and weak convergence to stochastic integrals II.}
\newblock{Econometric Theory,} 2000, 16: 643--666.\\[-6.5mm]
\bibitem{Dai2011}
Dai  H.
\newblock{Convergence in law to operator fractional Brownian motions.}
\newblock{J Theor Probab,} 2013, 26: 676--696.\\[-6.5mm]
\bibitem{David2004}
David  A.
\newblock{ L\'evy Processes and Stochastic Calculus.}
\newblock{ Cambridge: Cambridge University Press,} 2004.\\[-6.5mm]
\bibitem{dh1970}
de Haan  L.
\newblock{On regular variation and its application to the weak convergence of sample extremes.}
\newblock{Amsterdam: Math Centre, }1970.\\[-6.5mm]

\bibitem{DP2011}
Didier  G, Pipiras V.
\newblock{Integral representations and properties of operator fractional Brownian motions.}
\newblock{ Bernoulli, } 2011, 17: 1--33.\\[-6.5mm]
\bibitem{DP2011a}
Didier  G, Pipiras V.
\newblock{Exponents, symmetry groups and classification of operator fractional Brownian motions.}
\newblock{ J Theor Probab, } 2012, 25: 353--395.\\[-6.5mm]
\bibitem{EK1986}
Ethier  S,  Kurtz T.
\newblock{Markov Processes: Characterization and
Convergence.}
\newblock{ New York: John Wiley and Sons, } 1986.\\[-6.5mm]
\bibitem{HM1982}
Hudson W N, Mason  J D.
\newblock{ Operator-self-similar processes in a finite-dimensional space.}
\newblock{Trans  Amer Math  Soc,}  1982, 273: 281--297.\\[-6.5mm]
\bibitem{LL}
Lamperti  L.
\newblock{Semi-stable stochastic processes.}
\newblock{Trans  Amer Math Soc,} 1962, 104: 62-78.\\[-6.5mm]
\bibitem{LR1982}
Laha T L, Rohatgi V K.
\newblock{Operator self-similar processes in $\R^d$.}
\newblock{Stoch Process  Appl,} 1982, 12: 73--84.\\[-6.5mm]

\bibitem{MR2000}
Marinucci D, Robinson  P.
\newblock{Weak convergence of multivariate fractional processes.}
\newblock{Stoch  Process Appl,} 2000, 86: 103--120.\\[-6.5mm]
\bibitem{MX1999}
Mason J D , Xiao Y.
\newblock{Sample path properties of operator-self-similar Gaussian random fields.}
\newblock{Theory Probab Appl,} 2002, 46: 58--78.\\[-6.5mm]
\bibitem{MS2001}
Meerschaert  M M, Scheffler  H P.
\newblock{Limit Distributions for Sums of Independent Random Vectors: Heavy Tails in Theory and Practice.}
\newblock{ New York: John Wiley and Sons,} 2001.\\[-6.5mm]
\bibitem{ST1982}
Sato  K.
\newblock{Self-similar processes with independent increments.}
\newblock{Probab Th Rel Fields,} 1991, 89: 285--300.\\[-6.5mm]
\bibitem{ST94}
Samorodnitsky  G,  Taqqu M S.
\newblock{Stable Non-Gaussian Random Processes: Stochastic Models with Infinite Variance.}
\newblock{New York: Chapman and Hall, }1994.\\[-6.5mm]
\bibitem{S1993}
S\'{a}nchez de Naranjo  M V.
\newblock{Non-central limit theorems for non-linear functionals of $k$ Gaussian fields.}
\newblock{J Multivariate Anal,} 1993, 44: 227--255.\\[-6.5mm]
\bibitem{S1995}
S\'{a}nchez de Naranjo  M V.
\newblock{A central limit theorem for non-linear functionals of stationary Gaussian vector processes.}
\newblock{Stat Prob Lett,}  1993, 22: 223--230.\\[-6.5mm]
\bibitem{T1975}
Taqqu M S.
\newblock{ Weak convergence to fractional Brownian
motion and to the Rosenblatt process.}
\newblock{Z Wahrsch Verw
Gebiete,}  1975, 31: 287--302.\\[-6.5mm]
\bibitem{VW}
Vervaat W.
\newblock{Sample path properties of self-similar processes with stationary increments.}
\newblock{The Annals of Probability,} 1985, 13: 1--27.\\[-6.5mm]
}
\end{enumerate}
\end{document}